\documentclass{gtart}

\def\ifplaintex{\expandafter\ifx\csname documentclass\endcsname\relax}

\ifplaintex 
\else
\fi

\expandafter\ifx\csname beginpicture\endcsname\relax
\expandafter\ifx\csname documentclass\endcsname\relax
\input pictex \else
\input prepictex \input pictex \input postpictex \fi\fi

\def\gt{{\mathsurround=0pt\it $\cal G\mskip-2mu$eometry \&\ 
$\cal T\!\!$opology}}        

\def\gtp{{\mathsurround=0pt\it $\cal G\mskip-2mu$eometry \&\ 
$\cal T\!\!$opology $\cal P\!$ublications}}  

\def\lognumber#1{\def\thelognumber{#1}}
\def\volumenumber#1{\def\thevolumenumber{#1}}
\def\papernumber#1{\def\thepapernumber{#1}}
\def\volumeyear#1{\def\thevolumeyear{#1}}

\def\pagenumbers#1#2{\def\startpage{#1}\def\finishpage{#2}}
\def\published#1{\def\publishdate{#1}}
\def\proposed#1{\def\theproposer{#1}}
\def\seconded#1{\def\theseconders{#1}}
\def\received#1{\def\receiveddate{#1}}

\def\accepted#1{\def\accepteddate{#1}}
\def\asciititle#1{\def\theasciititle{#1}}

\def\coverauthors#1{\def\thecoverauthors{#1}}
\def\asciiauthors#1{\def\theasciiauthors{#1}}
\def\asciiaddress#1{\def\theasciiaddress{#1}}
\def\asciiemail#1{\def\theasciiemail{#1}}
\long\def\asciiabstract#1{\long\def\theasciiabstract{#1}}

\let\\\par\let\thelognumber\relax
\let\thevolumenumber\relax\let\thepapernumber\relax
\let\thevolumeyear\relax\let\thesamplenumber\relax\let\startpage\relax
\let\finishpage\relax\let\publishdate\relax\let\receiveddate\relax
\let\reviseddate\relax\let\accepteddate\relax\let\theasciititle\relax
\let\theasciiauthors\relax\let\theasciiaddress\relax
\let\theasciiabstract\relax
\let\theasciiemail\relax\let\theshortauthors\relax\let\theshorttitle\relax
\let\thecoverauthors\relax

\long\def\maketitlep{   

\count0=\startpage

\gt\hfill      
\beginpicture
\setcoordinatesystem units <0.33truein, 0.33truein> point at 2.2 0.9
\setplotsymbol ({$\cal G$})
\plotsymbolspacing=9truept
\circulararc 315 degrees from 0 1 center at 0 0
\setplotsymbol ({$\cal T$})
\circulararc 315 degrees from 1 -1 center at 1 0
\endpicture
\break
{\small\ifx\thesamplenumber\relax 
Volume \else Sample
\fi\thevolumenumber\ (\thevolumeyear)
\startpage--\finishpage\nl
Published: \publishdate}
\vglue 0.5truein plus 0.4fil minus 0.1truein

{\parskip=0pt\leftskip 0pt plus 1fil\def\\{\par\smallskip}{\ifplaintex\large
\else\Large\fi\bf\thetitle}\par\medskip}   

\vglue 0pt plus 0.1fil 

{\parskip=0pt\leftskip 0pt plus 1fil\def\\{\par}{\sc\theauthors}
\par\medskip}

\vglue 0pt plus 0.1fil 

{\small\parskip=0pt\let\newline\\
{\leftskip 0pt plus 1fil\def\\{\par}{\sl\theaddress}\par}
\expandafter\ifx\theemail\relax    
\relax\else\vglue 5pt plus 0.02fil minus 2pt\def\\{\stdspace{\rm 
and}\stdspace} 
\cl{Email:\stdspace\tt\theemail}\fi
\ifx\theurl\relax                  
\relax\else\vglue 5pt plus 0.02fil minus 2pt\def\\{\stdspace{\rm 
and}\stdspace}
\cl{URL:\stdspace\tt\theurl}\fi\par}

\vglue 7pt plus 0.3fil minus 3pt

{\bf Abstract}
\vglue 5pt plus 0.1fil minus 2pt

\theabstract

\vglue 7pt plus 0.3fil minus 3pt

{\bf AMS Classification numbers}\quad Primary:\quad \theprimaryclass

Secondary:\quad \thesecondaryclass

\vglue 5pt plus 0.3fil minus 2pt

{\bf Keywords}\quad \thekeywords

\vglue 10pt plus 0.5fil minus 5pt

{\small  Proposed: \theproposer\hfill Received: \receiveddate\nl
Seconded: \theseconders\hfill 
\ifx\reviseddate\relax                         
Accepted: \accepteddate                        
\else
Revised: \reviseddate                          
\fi}
\eject
}       

\let\maketitlepage\maketitlep
\let\maketitle\maketitlepage

\font\phead=cmsl9 scaled 950
\font=cmsl9 scaled 1050
\font\pnum=cmbx10 scaled 913
\font\lnum=cmbx10 
\font\pfoot=cmsl9 scaled 950
\font=cmsl9 scaled 1050
\ifplaintex
\headline{\vbox to 0pt{\vskip -4.5mm\line{\small\phead\ifnum
\count0=\startpage ISSN 1364-0380 (on line)
1465-3060 (printed) \hfill {\pnum\folio}\else\ifodd\count0\def\\{ }%
\ifx\theshorttitle\relax\thetitle\else\theshorttitle\fi\hfill{\pnum\folio}
\else\def\\{ and }{\pnum\folio}\hfill\ifx\theshortauthors\relax\theauthors
\else\theshortauthors\fi\fi\fi}\vss}}
\footline{\vbox to 0pt{\vglue 0mm\line{\small\pfoot\ifnum\count0=\startpage
\copyright\ \gtp\hfill\else
\gt, Volume \thevolumenumber\ (\thevolumeyear)\hfill\fi}\vss
}}
\else
\makeatletter
\def\@oddhead{{\small\lhead\ifnum\count0=\startpage ISSN 1364-0380 (on line)
1465-3060 (printed) \hfill {\lnum\number\count0}\else\ifodd\count0
\def\\{ }\ifx\theshorttitle\relax \thetitle \else\theshorttitle\fi\hfill
{\lnum\number\count0}\else\def\\{ and }{\lnum\number\count0}
\hfill\ifx\theshortauthors\relax 
\theauthors\else\theshortauthors\fi\fi\fi}}\def\@evenhead{\@oddhead}
\def\@oddfoot{\small\lfoot\ifnum\count0=\startpage\copyright\ \gtp\hfill\else
\gt, Volume \thevolumenumber\ (\thevolumeyear)\hfill\fi}
\def\@evenfoot{\@oddfoot}
\makeatother
\fi

\newwrite\gtoutfile
\long\gdef\makeheadfile{  
{\def\\{, }\def\s{ }
\immediate\openout\gtoutfile head.xxx
\immediate\write\gtoutfile{Proxy-for: \ifx\theasciiauthors\relax
\theauthors\else\theasciiauthors\fi\s<\ifx\theasciiemail\relax\theemail\else\theasciiemail\fi>}
\immediate\write\gtoutfile{\noexpand\\}
\immediate\write\gtoutfile{Authors: \ifx\theasciiauthors\relax
\theauthors\else\theasciiauthors\fi}
{\def\\{ }\immediate\write\gtoutfile{Title: \ifx\theasciititle\relax
\thetitle\else\theasciititle\fi}}
\immediate\write\gtoutfile{Subj-class: GT or SG or MG etc}
\immediate\write\gtoutfile{MSC-class: \theprimaryclass\ifx\thesecondaryclass\relax\else, \thesecondaryclass\fi}
\immediate\write\gtoutfile{Journal-ref: Geom. Topol. \thevolumenumber
(\thevolumeyear) \startpage-\finishpage}
\immediate\write\gtoutfile{Comments: Published by Geometry and Topology at}
\immediate\write\gtoutfile{\s\s http://www.maths.warwick.ac.uk/gt/GTVol\thevolumenumber/paper\thepapernumber.abs.html}
\immediate\write\gtoutfile{\noexpand\\}
\immediate\write\gtoutfile{}
\ifx\theasciiabstract\relax
\immediate\write\gtoutfile{\theabstract}\else
\immediate\write\gtoutfile{\theasciiabstract}\fi
\immediate\write\gtoutfile{}
\immediate\write\gtoutfile{\noexpand\\}
\immediate\write\gtoutfile{}
\immediate\closeout\gtoutfile}}  

\def\maketitlepage{\maketitlep\makeheadfile}
\let\maketitle\maketitlepage

\lognumber{272}
\volumenumber{7}\papernumber{30}\volumeyear{2003}
\pagenumbers{1055}{1073}
\received{4 September 2002}
\published{29 December 2003}
\accepted{23 December 2003}
\proposed{Tomasz Mrowka}
\seconded{Yasha Eliashberg, John Morgan}

\usepackage{amssymb, amsmath}

\newtheorem{thm}{Theorem}[section]

\newtheorem{lem}[thm]{Lemma}
\newtheorem{prop}[thm]{Proposition}

\theoremstyle{definition}
\newtheorem{defn}[thm]{Definition}

\newtheorem{rem}[thm]{Remark}

\numberwithin{equation}{section}

\newcommand{\al}{\alpha}

\newcommand{\Ga}{\Gamma}

\newcommand{\ep}{\epsilon}

\newcommand{\La}{\Lambda}
\newcommand{\la}{\lambda}
\newcommand{\Om}{\Omega}
\newcommand{\om}{\omega}

\newcommand{\Si}{\Sigma}

\newcommand{\ze}{\zeta}
\newcommand{\bfz}{{\mathbb {Z}}}

\newcommand{\x}{\times}
\newcommand{\s}{\mathbf s}

\renewcommand{\t}{\mathbf t}

\renewcommand{\L}{\mathbf\La}
\newcommand{\D}{\mathbf D}

\newcommand{\C}{\mathbb C}
\newcommand{\Z}{\mathbb Z}

\newcommand{\R}{\mathbb R}
\newcommand{\bfr}{\mathbb R}

\newcommand{\CP}{{\mathbb C}{\mathbb P}}

\newcommand{\del}{\partial}

\newcommand{\lra}{\longrightarrow}

\renewcommand{\rk}{\operatorname{\mathrm rk}}

\newcommand{\PD}{\operatorname{\mathrm PD}}

\begin{document}
\title{An infinite family of tight, not
semi--fillable\\contact three--manifolds}
\asciititle{An infinite family of tight, not
semi-fillable contact three-manifolds}

\author{Paolo Lisca\\Andr\'{a}s I Stipsicz}
\asciiauthors{Paolo Lisca\\Andras I Stipsicz}
\coverauthors{Paolo Lisca\\Andr\noexpand\'{a}s I Stipsicz}
\address{Dipartimento di Matematica,
Universit\`a di Pisa\\I-56127 Pisa, ITALY} 
\email{lisca@dm.unipi.it}

\secondaddress{R\'enyi Institute of Mathematics,
Hungarian Academy of Sciences\\H-1053 Budapest, 
Re\'altanoda utca 13--15, Hungary}
\secondemail{stipsicz@math-inst.hu}

\asciiemail{lisca@dm.unipi.it, stipsicz@math-inst.hu}
\asciiaddress{Dipartimento di Matematica,
Universita di Pisa\\I-56127 Pisa, ITALY\\and\\Renyi 
Institute of Mathematics,
Hungarian Academy of Sciences\\H-1053 Budapest, 
Realtanoda utca 13--15, Hungary}

\begin{abstract}
  We prove that an infinite family of virtually overtwisted tight
  contact structures discovered by Honda on certain circle bundles
  over surfaces admit no symplectic semi--fillings. The argument uses
  results of Mrowka, Ozsv\'ath and Yu on the translation--invariant
  solutions to the Seiberg--Witten equations on cylinders and the
  non--triviality of the Kronheimer--Mrowka monopole invariants of 
  symplectic fillings.
\end{abstract}
\asciiabstract{%
  We prove that an infinite family of virtually overtwisted tight
  contact structures discovered by Honda on certain circle bundles
  over surfaces admit no symplectic semi-fillings. The argument uses
  results of Mrowka, Ozsvath and Yu on the translation-invariant
  solutions to the Seiberg-Witten equations on cylinders and the
  non-triviality of the Kronheimer-Mrowka monopole invariants of 
  symplectic fillings.}
\primaryclass{57R57}
\secondaryclass{57R17}
\keywords{Tight, fillable, contact structures}

\maketitlepage

\section{Introduction}\label{s:intro}
Let $Y$ be a closed, oriented three--manifold. A \emph{positive,
  coorientable contact structure} on $Y$ is the kernel
$\xi=\ker\al\subset TY$ of a one--form $\alpha\in\Om^1(Y)$ such that
$\alpha\wedge d\alpha$ is a positive volume form on $Y$. In this paper
we only consider positive, coorientable contact structures, so we call
them simply `contact structures'. For an introduction to contact
structures the reader is referred to \cite[Chapter~8]{Ae}
and~\cite{E}.

There are two kinds of contact structures $\xi$ on $Y$. If there
exists an embedded disk $D\subset Y$ tangent to $\xi$ along its
boundary, $\xi$ is called \emph{overtwisted}, otherwise it is said to
be \emph{tight}. The isotopy classification of overtwisted contact
structures coincides with their homotopy classification as tangent
two--plane fields~\cite{Eli}. Tight contact structures are much more
difficult to classify, and capture subtle information about the
underlying three--manifold. 

A contact three--manifold $(Y,\xi)$ is \emph{symplectically fillable},
or simply \emph{fillable}, if there exists a compact symplectic
four--manifold $(W, \omega)$ such that (i) $\partial W=Y$ as oriented
manifolds (here $W$ is oriented by $\om\wedge\om$) and (ii)
$\omega\vert_{\xi}\not=0$ at every point of $Y$. $(Y,\xi)$ is
symplectically \emph{semi}--fillable if there exists a fillable
contact manifold $(N, \eta )$ such that $Y\subset N$ and
$\eta|_Y=\xi$. Semi--fillable contact structures are
tight~\cite{Eli1}. The converse is known to be false by work of Etnyre
and Honda, who recently found two examples of tight but not
semi--fillable contact three--manifolds~\cite{EH}. This discovery
naturally led to a search for such examples, in the hope that they
would tell us something about the difference between tight and
fillable contact structures. By a result announced by
E~Giroux~\cite{Gi3}, isotopy classes of contact structures on a closed
three--manifold are in one--to--one correspondence with ``stable''
isotopy classes of open book decompositions. When the monodromy of the
open book decomposition is positive, the corresponding contact
structure is fillable. Therefore, it would be very interesting to know
examples of monodromies associated with tight but not fillable contact
structures.

In this paper we prove that infinitely many tight contact circle
bundles over surfaces are not semi--fillable. Let $\Si_g$ be a closed,
oriented surface of genus $g\geq 1$. Denote by $Y_{g,n}$ the total
space of an oriented $S^1$--bundle over $\Si_g$ with Euler number $n$.
Honda gave a complete classification of the tight contact structures
on $Y_{g,n}$~\cite{H2}. The three--manifolds $Y_{g,n}$ carry
infinitely many tight contact structures up to diffeomorphism.  The
classification required a special effort for two tight contact
structures $\xi_0$ and $\xi_1$, for $n$, $g$ satisfying $n\geq 2g$
(see Definition~\ref{d:contact} below). Honda conjectured that none of
these contact structures are symplectically semi--fillable. Our main
result establishes Honda's conjecture in infinitely many cases:

\begin{thm}\label{t:main}
Suppose that $d(d+1)\leq 2g\leq n\leq (d+1)^2+3$ for some positive
integer $d$. Then, the tight contact structures $\xi_0$ and $\xi_1$ on
$Y_{g,n}$ are not symplectically semi--fillable.
\end{thm}

\begin{rem}
  Let $\xi$ be a tight contact structure on a three--manifold $Y$.  If
  the pull--back of $\xi$ to the universal cover of $Y$ is tight, then
  $\xi$ is called~\emph{universally tight}. Honda showed in~\cite{H2}
  that the contact structures $\xi_i$ become overtwisted when
  pulled--back to a finite cover, i.e.\ they are \emph{virtually
    overtwisted}. Thus, the question whether every universally tight
  contact structure is symplectically fillable is untouched by
  Theorem~\ref{t:main}.
\end{rem}

The proof of Theorem~\ref{t:main} is similar in spirit to the argument
used by the first author to prove that certain oriented
three--manifolds with positive scalar curvature metrics do not carry
semi--fillable contact structures~\cite{PLpos, Paolo}. But the fact
that the three--manifolds $Y_{g,n}$ do not admit positive scalar
curvature metrics made necessary a modification of the analytical as
well as the topological parts of the original argument.

We first show that if $W$ is a semi--filling of $(Y_{g,n}, \xi _i)$,
then $\partial W$ is connected, $b_2^+(W)=0$ and the homomorphism $H^2
(W; \bfr )\to H^2 (\partial W;\bfr )$ induced by the inclusion
$\partial W \subset W$ is the zero map (see
Proposition~\ref{p:filling}). To do this, we start by identifying the
$Spin^c$ structures $\t_{\xi_i}$ induced by the contact structures
$\xi_i$. Then, using results of Mrowka, Ozsv\'ath and Yu~\cite{MOY},
we establish properties of the Seiberg--Witten moduli spaces for the
$Spin^c$ structures $\t_{\xi_i}$ which are sufficient to apply the
argument used in the positive scalar curvature case. Such an argument
relies on the non--triviality of the Kronheimer--Mrowka monopole
invariants of a symplectic filling~\cite{KrMr}.

Then, under some restrictions on $g$ and $n$, we construct smooth,
oriented four--manifolds $Z$ with boundary orientation--reversing
diffeomorphic to $Y_{g,n}$, with the property that if $W$ were a
symplectic filling of $(Y_{g,n}, \xi_i)$, the closed four--manifold
$W\cup _{Y_{g,n}}Z$ would be negative definite and have
non--diagonalizable intersection form. On the other hand, by
Donaldson's celebrated theorem~\cite{Do1, Do2}, such a closed four--manifold
cannot exist. Therefore, $(Y_{g,n}, \xi_i)$ does not have symplectic
fillings.

The paper is organized as follows. In Section~\ref{s:definition} we
define, following~\cite{H2}, the contact structures $\xi_0$ and
$\xi_1$. In Section~\ref{s:spinccalc} we determine the $Spin ^c$
structures $\t_{\xi_i}$, and in Section~\ref{s:gauge} we prove
Theorem~\ref{t:main}.

{\bf Acknowledgements}\qua The first author is grateful to Peter Ozsv\'ath
for enlighting e--mail correspondence. The first author was partially
supported by MURST, and he is a member of EDGE, Research Training
Network HPRN-CT-2000-00101, supported by The European Human Potential
Programme. The second author was partially supported by Bolyai
\"Oszt\"ond{\'\i}j and OTKA T34885.

\section{Definition of the contact structures}
\label{s:definition}

In this section we describe in detail the construction of the contact
structures $\xi _i$. For the sake of the exposition, we start by
recalling some basic facts regarding convex surfaces and Legendrian
knots in contact three--manifolds.

\subsection*{Basic properties of contact structures}
Let $Y$ be a closed, oriented three--manifold and let $\xi$ be a contact
structure on $Y$. 

\begin{defn}
An embedded surface $\Si\subset Y$ is \emph{convex} if there exists a
vector field $V$ on $Y$ such that (i) $V$ is transverse to $\Si$ and
(ii) $V$ is a \emph{contact vector field}, i.e.\ $\xi$ is invariant
under the flow generated by $V$. The \emph{dividing set} is
\[
\Ga_\Si(V)=\{p\in\Si\ |\ V(p)\in\xi\}\subset\Si.
\]
\end{defn}

The following facts are proved in the seminal paper by
E~Giroux~\cite{Gi1}:
\begin{enumerate}
\item
Let $\Si\subset (Y,\xi)$ be an embedded surface. Then, $\Si$ can be
$C^\infty$--perturbed to a convex surface.
\item
Let $\Si\subset (Y,\xi)$ be a convex surface and $V$ a contact vector
field transverse to $\Si$. Then, (i) the isotopy class $\Ga_\Si$ of
$\Ga_\Si(V)$ does not depend on the choice of $V$ and (ii) the germ of
$\xi$ around $\Si$ is determined by $\Ga_\Si$.
\end{enumerate}

In the case of a convex torus $T\subset (Y,\xi)$, the set $\Ga_T$
consists of an even number of disjoint simple closed curves. The germ
of $\xi$ around $T$ is determined by the number of connected
components of $\Ga_T$ -- the \emph{dividing curves} -- together with a
(possibly infinite) rational number representing their slope with
respect to an identification $T\cong\R^2/\Z^2$.  Note that the slope
depends on the choice of the identification. If $T$ is the boundary of
a neighborhood of a knot $k\subset Y$, then by identifying the
meridian with one copy of $\bfr /\bfz $ in the above identification,
the slope $p/q$, regarded as a vector $\binom{q}{p}$, is determined up
to the action of the group 
\[
\{\left(\begin{matrix} 
1 & m \\ 
0 & 1
\end{matrix}\right) \mid m \in \bfz \}.
\] 
A knot $k\subset (Y,\xi)$ is \emph{Legendrian} if $k$ is everywhere
tangent to $\xi$. The framing of $k$ naturally induced by $\xi$ is called
the \emph{contact framing}. A Legendrian knot $k$ has a basis $\{U_\al\}$
of neighborhoods with convex boundaries. The dividing set of each boundary
$\del U_\al$ consists of two dividing curves having the same slope
independent of $\al$. Any one of those neighborhoods of $k$ is called a
\emph{standard convex neighborhood} of $k$. The meridian and the contact
framing of a Legendrian knot $k \subset (Y, \xi )$ provides an
identification of the convex boundary $T$ of a neighborhood of $k$ with
$\bfr ^2/\bfz ^2$; easy computation shows that with this identification
the slope of the dividing curves is $\infty$.

Let $\Si_g$ be a closed, oriented surface of genus $g\geq 1$, and let
$\pi\co Y_{g,n}\to\Si_g$ be an oriented circle bundle over $\Si_g$
with Euler number $n$. Let $\xi$ be a contact structure on $Y_{g,n}$
such that a fiber $f=\pi^{-1}(s)\subset Y_{g,n}$ ($s\in\Si_g$) is
Legendrian. We say that $f$ has \emph{twisting number $-1$} if the
contact framing of $f$ is `$-1$' with respect to the framing
determined by the fibration $\pi$.  A contact structure on $Y_{g,n}$
is called \emph{horizontal} if it is isotopic to a contact structure
transverse to the fibers of $\pi$.

\subsection*{Definition of the contact structures $\xi_i$}

The following lemma is probably well--known to the experts. A proof
can be found e.g.~in~\cite[\S 1.D]{Gi2}. We include it here to make
our exposition more self--contained, and for later reference.

\begin{lem}\label{l:horizontal}
The circle bundle $\pi\co Y_{g,2g-2}\to\Si_g$ carries a horizontal
contact structure $\ze$ such that all the fibers of $\pi$ are
Legendrian and have twisting number $-1$.
\end{lem}

\begin{proof}
Think of $Y_{g,2g-2}$ as the manifold of the oriented lines tangent to
$\Si_g$. Then, the fiber $\pi^{-1}(s)\subset Y_{g,2g-2}$ ($s\in\Si_g$)
consists of all the oriented lines $l$ tangent to $\Si_g$ at the point
$s$. The contact two--plane $\ze(l)$ at the point $l$ is, by definition, the
preimage of $l\subset T_s \Si_g$ under the differential of $\pi$. It is a
classical fact that $\ze$ is a contact structure. It follows directly from the
definition that every fiber of $\pi$ is Legendrian and has twisting number
$-1$. To see that $\ze$ is horizontal, let $V$ be a vector field on
$Y_{g,2g-2}$ tangent to $\ze$, transverse to the fibers and such that, for
every $l\in Y_{g,2g-2}$ the projection $d\pi(V(l))\in l\subset T_{\pi(l)}
Y_{g,2g-2}$ defines the orientation on $l$. Let $\la$ be a one--form defining
$\ze$, and $T$ a nonvanishing vector field on $Y_{g,2g-2}$ tangent to the
fibers. The fact that $d\la$ does not vanish on the contact planes is
equivalent to the fact that the Lie derivative of $\la$ in the direction of
$V$ is nowhere vanishing when evaluated on $T$, because ${\mathfrak
{L}}_V(\la)(T)=d\la(V,T)$. Thus, following the flow of $V$ the contact
structure $\ze$ can be isotoped to a transverse contact structure.
\end{proof}
Suppose that $n\geq 2g$. Let $\ze$ be the contact structure given by
Lemma~\ref{l:horizontal} (according to~\cite[\S 5]{H2} $\ze$ could be any
horizontal contact structure on $Y_{g,2g-2}$ such that a fiber $f$ of
the projection $\pi$ is Legendrian with twisting number $-1$). Let
$U\subset Y_{g,2g-2}$ be a standard convex neighborhood of $f$. The
fibration induces a trivialization $U\cong S^1\x D^2$. Remove $U$ from
$Y_{g,2g-2}$ and reglue it using the diffeomorphism $\varphi_A\co\del
U\to -\del (Y_{g,2g-2}\setminus U)$ determined, via the above
trivialization, by the matrix
\[
A=
\begin{pmatrix}
  -1 & 0 \\
  p+1 & -1             
\end{pmatrix},
\]
where $p=n-2g+1$. The map $\pi$ extends to the resulting
three--manifold yielding the bundle $Y_{g,n}$, and we are going to
show that $\ze$ extends as well. The germ of $\ze$ around
$\del(Y_{g,2g-2}\setminus U)$ is determined by the slope of any
dividing curve $C\subset\del(Y_{g,2g-2}\setminus U)$. We are going to
extend $\ze$ to $U$ as a tight contact structure $\eta$ having convex
boundary and two dividing curves isotopic to
$\varphi_A^{-1}(C)\subset\del U$. Since the fiber $f$ has twisting
number $-1$ with respect to $\ze$, the slope of $C$ is $-1$ with
respect to the trivialization used to define the gluing
map. Therefore, the slope of $\varphi_A^{-1}(C)$ is $p$. Applying the
self--diffeomorphism of $S^1\x D^2$ given by the matrix
$\left(\begin{smallmatrix} 1 & -1 \\ 0 & \phantom{-}1
\end{smallmatrix}
\right)$, we may assume that the boundary slope of $(U,\eta)$ be
$-\frac p{p-1}$.

Summarizing, assuming the existence of $\eta$ we have constructed
a contact three--manifold $(Y_{g,n},\xi)$ of the form
\begin{equation}\label{e:contactgluing}
(Y_{g,n},\xi)=(Y_{g,2g-2}\setminus U,\ze)\cup_{\varphi_A}
(U,\eta).
\end{equation}
The following result says that we have two possible choices for $\eta$
when $n > 2g$, and one when $n=2g$.

\begin{thm}\label{t:s1xd2}{\rm\cite{H1}}\qua 
Let $p$ be a positive integer. Up to an isotopy keeping the boundary
fixed, there are at most two tight contact structures on $S^1\times D^2$
with convex boundary and two dividing curves with slope $-\frac
p{p-1}$. More precisely, when $p>1$ there are exactly two such contact
structures. When $p=1$, i.e.\ when the boundary slope is infinite,
there is only one. \qed
\end{thm}

We need to be more specific about the contact structures appearing in
the statement of Theorem~\ref{t:s1xd2}. The next lemma will be used in
the following section as well.

\begin{lem}\label{l:torus-section}
Let $p$ be a positive integer and $\eta$ a tight contact structure on
$S^1\times D^2$ with convex boundary and two dividing curves with
slope $-\frac{p}{p-1}$. Then, $\eta$ can be isotoped keeping the
boundary convex until there exists a section $v$ of $\eta$ such that:
\begin{enumerate}
\item
At every boundary point $v$ is nonvanishing and tangent to 
the circles $S^1\x\{x\}$, $x\in\del D^2$;
\item
The zero locus of $v$ is a smooth curve homologous to 
\begin{equation}\label{e:zerolocus}
\pm (p-1)[S^1\x\{(0,0)\}].
\end{equation}
\end{enumerate}
Moreover, each sign in formula~\eqref{e:zerolocus} above can be
realized by some contact structure $\eta$.
\end{lem}

\begin{proof}
Using Giroux's Flexibility Theorem (see~\cite{Gi1} and~\cite[\S 3.1.4]{H1}) 
we may isotope $\eta$ keeping the boundary convex until $\eta$
is tangent to the circles $S^1\x\{x\}$, $x\in\del D^2$ at each
boundary point.

By \cite[Proposition~4.15]{H1} there exists a decomposition 
\[
S^1\x D^2=N\cup (S^1\x D^2\setminus N), 
\]
where $N$ is a standard convex neighborhood of a Legendrian knot
isotopic to the core circle of $S^1\x D^2$. Thus, $N\cong S^1\x D^2$
with coordinates $(z,(x,y))$ and
\[
\eta|_N=\ker\left(\sin(2\pi z)dx -\cos(2\pi z)dy\right).
\]
Moreover, there is a diffeomorphism 
\[
\varphi\co S^1\x D^2\setminus N\cong T^2\x [0,1]
\]
and $(T^2\x [0,1],\varphi_*\eta)$ is a \emph{basic slice}
(see~\cite[\S 4.3]{H1}) with convex boundary components $T^2\x\{0\}$
and $T^2\x\{1\}$ of slopes $-1$ and $-\frac p{p-1}$
respectively. Without loss of generality we may also assume that
\[
\varphi|_{S^1\x\del D^2}\co S^1\x\del D^2\to T^2\x\{1\}
\]
is the obvious identification.

According to \cite[Lemma~4.6 and Proposition~4.7]{H1} the Euler class
of a basic slice with boundary slopes $0$ and $-1$, relative to a
section which is nowhere zero at the boundary and tangent to it, is
equal to 
\[
\pm (0,1)\in H_1(T^2\x [0,1];\Z)\cong H_1(T^2;\Z)\cong\Z^2.
\]
Moreover, each sign is realized by a unique (up to isotopy) basic
slice. Applying the diffeomorphism 
\[
T^2\times [0,1]\lra T^2 \times [0,1]
\]
given by 
\[
\left(\begin{matrix} 
\phantom{-}1 & 2-p \\ 
-1 & p-1
\end{matrix}\right)
\] 
we obtain a basic slice with boundary slopes
$-1$ and $-\frac p{p-1}$ which is diffeomorphic to $(S^1\x
D^2\setminus N,\eta)$, together with a section $v$ of $\eta$ which is
nowhere zero at the boundary, tangent to it, and with zero locus a
smooth curve homologous to $\pm (2-p,p-1)$. The section $v$ can be
assumed to coincide with the vector field $\frac{\del}{\del z}$ on
$\del N$ in the above coordinates $(z,(x,y))$. Since $\frac{\del}{\del
z}$ is a nowhere zero section of $\eta$ on $N$, this implies that $v$
extends as a section of $\eta$ on $S^1\x D^2$ with the stated
properties. By choosing the appropriate basic slice one can construct
$\eta$ with either choice of sign in Formula~\eqref{e:zerolocus}.
\end{proof}

\begin{defn}\label{d:contact}
  Let $\eta_0$ (respectevely $\eta_1$) be a tight contact structure on
  $S^1\x D^2$ as in the conclusion of Lemma~\ref{l:torus-section},
  satisfying condition~\eqref{e:zerolocus} with the positive
  (respectively negative) sign. Let $n\geq 2g$, and define $\xi_0$
  (respectively $\xi_1$) to be the contact structure on $Y_{g,n}$,
  $n\geq 2g$, given by~\eqref{e:contactgluing} with $\eta$ replaced by
  $\eta_0$ (respectively $\eta_1$).
\end{defn}

\begin{rem}\label{r:xi}
By Theorem~\ref{t:s1xd2} the contact structures $\xi_0$ and $\xi_1$
of Definition~\ref{d:contact} are isotopic when $n=2g$. By the
classification from~\cite{H2}, $\xi_0$ is not isotopic to $\xi_1$ when
$n>2g$. In fact, $\xi_0$ and $\xi_1$ are not even homotopic (see
Remark~\ref{r:homotopy} below). 
\end{rem}

\section{Calculations of $Spin^c$ structures}
\label{s:spinccalc}
The goal of this section is to determine the $Spin^c$ structures
$\t_{\xi_i}$ induced by the contact structures $\xi_i$ of
Definition~\ref{d:contact}. We begin with a short review about $Spin$
and $Spin^c$ structures in general. Then, we study $Spin$ and $Spin^c$
structures on disk and circle bundles over surfaces. The section ends
with the calculation of $\t_{\xi_0}$ and $\t_{\xi_1}$.

\subsection*{Generalities on $Spin$ and $Spin^c$ structures}
\begin{sloppypar}
Let $X$ be a smooth, oriented $n$--dimensional manifold, $n\geq
3$. The structure group of its tangent frame bundle $P_X$ can be
reduced to $SO(n)$ by e.g.\ introducing a Riemannian metric on $X$. A
\emph{$Spin$ structure} on $X$ is a principal $Spin(n)$--bundle
$P_{Spin(n)}\to X$ such that $P_X$ is isomorphic to the associated
bundle
\[
P_{Spin(n)}\x_\rho SO(n),
\]
where 
\[
\rho\co Spin(n)\to SO(n)
\]
is the universal covering map. A $Spin$ structure on $X$ exists if and
only if the second Stiefel--Whitney class $w_2(X)$ vanishes. In this
case, the set of $Spin$ structures is a principal homogeneous space on
$H^1(X;\Z/2\Z)$.
\end{sloppypar}

The quotient of $Spin(n)\x S^1$ modulo the subgroup 
\[
\{\pm (1,1)\}\cong\Z/2\Z
\]
is, by definition, the group $Spin^c(n)$.  There
are two canonical surjective homomorphisms 
\[
\rho_1\co Spin^c(n)\to Spin(n)/\{\pm 1\}= SO(n),\qquad
\rho_2\co Spin^c(n)\to S^1/\{\pm 1\}= S^1. 
\]
A \emph{$Spin^c$ structure} on $X$ is a principal $Spin^c(n)$--bundle
$P_{Spin^c(n)}$ such that 
\[
P_X\cong P_{Spin^c(n)}\x_{\rho_1} SO(n).
\]
Let $Spin^c(X)$ denote the (possibly empty) set of $Spin^c$ structures
on $X$. An element 
\[
\s=P_{Spin^c(n)}\in Spin^c(X)
\]
naturally induces a principal $S^1$--bundle $P_{Spin^c(n)}\x_{\rho_2}
S^1$. Let $c_1(\s)\in H^2(X;\Z)$ be the first Chern class of the
corresponding complex line bundle. A manifold $X$ admits a $Spin^c$
structure if and only if $w_2(X)$ has an integral lift, and in fact
the set 
\[
\{c_1(\s)\ |\ \s\in Spin^c(X)\}\subset H^2(X;\Z)
\]
is the preimage of $w_2(X)$ under the natural map 
\[
H^2(X;\Z)\to H^2(X;\Z/2\Z). 
\]
Moreover, $Spin^c(X)$ is a principal homogeneous space
on $H^2(X;\Z)$, and $c_1(\s+\al)=c_1(\s)+2\al$ for every 
\[
(\s,\al)\in Spin^c(X)\x H^2(X;\Z).
\]
The group $Spin(n)$ naturally embeds into $Spin^c(n)$, so a
$Spin$ structure induces a $Spin^c$ structure. Moreover, since
\[
Spin(n)=\ker\rho_2\subset Spin^c(n), 
\]
a $Spin^c$ structure $\s$ is induced by a $Spin$ structure if and only
if $c_1(\s)=0$.

Since 
\[
\rho_1^{-1}(\{1\}\x SO(n))=Spin^c(n)\subset Spin^c(n+1),
\]
if $\dim Y=n$ and $Y =\del X$, there is a restriction map
$Spin^c(X)\to Spin^c(Y)$. Clearly, this map sends $Spin$
structures to $Spin$ structures.

An oriented two--plane field $\xi$ (and so a contact structure) on a
closed, oriented three--manifold $Y$ reduces the structure group of
$TY$ to $U(1)\subset SO(3)$. Since the inclusion $U(1)\subset SO(3)$
admits a canonical lift to $U(2)=Spin^c(3)$, there is a $Spin^c$
structure $\t_\xi\in Spin^c(Y)$ canonically associated to $\xi$. The
$Spin^c$ structure $\t_\xi$ depends only on the homotopy class of
$\xi$ as an oriented tangent two--plane field.

\subsection*{Disk bundles}

Let $\Si_g$ be a closed, oriented surface of genus $g\geq 1$. Let
$\pi\co D_{g,n}\to\Si_g$ be an oriented $2$--disk bundle over $\Si_g$
with Euler number $n$. 

By e.g.\ fixing a metric on $D_{g,n}$ one sees that the tangent bundle
of $D_{g,n}$ is isomorphic to the direct sum of the pull--back of
$T\Si_g$ and the vertical tangent bundle, which is isomorphic to the
pull--back of the real oriented two--plane bundle $E_{g,n}\to\Si_g$
with Euler number $n$. In short, we have
\begin{equation}\label{e:split}
TD_{g,n}\cong \pi^*(T\Si_g\oplus E_{g,n}).
\end{equation}
Therefore, the structure group of $TD_{g,n}$ can be reduced to
$U(2)\subset SO(4)$, which admits a natural lift 
\[
A\mapsto (\left(
\begin{matrix}
\det A & 0 \\ 
\phantom{aa}0\phantom{aa} & 1
\end{matrix}\right) ,A)
\]
$$
Spin^c(4)=\{ (A,B)\in U(2)\times U(2) \mid \det A=\det B\}.\leqno{\hbox{to}}
$$
Denote by $\s_0$ the induced $Spin^c$ structure on $D_{g,n}$. The
orientation on $D_{g,n}$ determines an isomorphism
$H^2(D_{g,n};\Z)\cong\Z$, so the set 
\[
Spin^c(D_{g,n})=\s_0 + H^2(D_{g,n};\Z) 
\]
can be canonically identified with the integers. We
denote by 
\[
\s_e=\s_0+e\in Spin^c(D_{g,n})
\]
the element corresponding to the integer $e\in\Z\cong
H^2(D_{g,n};\Z)$.

\begin{lem}\label{l:spin4}
\noindent{\bf (a)}\qua
If $n$ is odd, $D_{g,n}$ admits no $Spin$ structure. 

\noindent{\bf (b)}\qua If $n$ is even, $D_{g,n}$ carries $Spin$
structures. Every $Spin$ structure on $D_{g,n}$ induces the $Spin^c$
structure $\s_{g-\frac n2-1}$.
\end{lem}

\begin{proof}
In view of~\eqref{e:split} we have $c_1(\s_0)=2-2g+n$, hence
\[
c_1(\s_e)=c_1(\s_0)+2e=2(1-g+e)+n. 
\]
Since each $c_1(\s_e)$ reduces modulo $2$ to $w_2(D_{g,n})$, $D_{g,n}$
admits a $Spin$ structure if and only if $n$ is even. Solving the
equation $c_1(\s_e)=0$ for $e$ yields the statement.
\end{proof}

\subsection*{Circle bundles}\label{ss:circle}

Consider $Y_{g,n}=\del D_{g,n}$. We have 
\[
H_1(Y_{g,n};\Z)\cong H^2(Y_{g,n}; \Z) \cong \Z
^{2g}\oplus \Z/n\Z, 
\]
where the summand $\Z/n\Z$ is generated by the
Poincar\'e dual $F$ of the class of a fiber of the projection $\pi\co
Y_{g,n}\to\Si_g$. Each $Spin^c$ structure $\s_e\in Spin^c(D_{g,n})$
determines by restriction a $Spin^c$ structure $\t_e\in
Spin^c(Y_{g,n})$. We have 
\[
\t_e=\t_0+eF,\quad e\in\Z. 
\]
Since $nF=0$, we see that $\t_{e+n}=\t_e$ for every $e$. Therefore,
$\t_0,\ldots,\t_{n-1}$ is a complete list of \emph{torsion $Spin^c$
structures} on $Y_{g,n}$, i.e.\ $Spin^c$ structures on $Y_{g,n}$ with
torsion first Chern class.  Notice that for $n$ even different
$Spin^c$ structures might have coinciding first Chern classes; for $n$
odd, $c_1(\t _i)$ determines $\t _i$.

\begin{rem}\label{r:transverse}
The pull--back of $E_{g,n}$ is trivial when restricted to the
complement of the zero section, therefore we have
\[
TY_{g,n}\cong\underline{\R}\oplus\pi^*(T\Si_g)\subset
TD_{g,n}|_{Y_{g,n}}\cong\underline{\C}\oplus\pi^*(T\Si_g),
\]
where $\underline{\R}$ and $\underline\C$ are, respectively, the
trivial real and complex line bundles. This shows that 
\[
\t_0=\s_0|_{Y_{g,n}}=\t_\ze,
\]
where $\ze\subset TY_{g,n}$ is any oriented tangent two--plane field
transverse to the fibers of $\pi\co Y_{g,n}\to\Si_g$.
\end{rem}

\begin{lem}\label{l:spin3}
\noindent{\bf (a)}\qua If $n$ is odd, $\t_{g-1}$ is the only torsion
$Spin^c$ structure on $Y_{g,n}$ induced by a $Spin$ structure.

\noindent{\bf (b)}\qua If $n$ is even, $\t_{g-1}$ and $\t_{g+\frac n2-1}$
are the only torsion $Spin^c$ structures on $Y_{g,n}$ induced by a
$Spin$ structure.
\end{lem}

\begin{proof}
Since $c_1(\s_0)$ restricts to $H^2(Y_{g,n};\Z)$ as $(2-2g)F$, we have
\[
c_1(\t_e)=2(1-g+e)F. 
\]
Solving the equation $c_1(\t_e)=0$ yields the statement.
\end{proof}

\begin{rem}\label{r:ext}
The $Spin^c$ structure $\t_{g-1}$ on $Y_{g,n}$ is (by definition) the
restriction of a $Spin^c$ structure on $D_{g,n}$. Although $\t_{g-1}$
is induced by a $Spin$ structure on $Y_{g,n}$, by Lemma~\ref{l:spin4}
$\t_{g-1}$ does not extend as a $Spin$ structure to $D_{g,n}$ when
$n>0$. On the other hand, when $n$ is even the $Spin^c$ structure
$\t_{g+\frac n2-1}$ is induced by a $Spin$ structure on $Y_{g,n}$
which is the restriction of a $Spin$ structure on $D_{g,n}$.
\end{rem}

\subsection*{Calculations}

Let $Y$ be a closed, oriented three--manifold and let $\xi\subset TY$ be an
oriented tangent two--plane field. The $Spin^c$ structure $\t_\xi\in
Spin^c(Y)$ determined by $\xi$ can be also defined as follows~\cite{G}.  Using
a trivialization of $TY$, the oriented two--plane bundle $\xi$ can be realized
as the pull--back of the tangent bundle to the two--sphere $S^2$ under a
smooth map $Y\to S^2$. This implies, in particular, that the Euler class of
$\xi$ is always even. Therefore $\xi$ has a section $v$ which vanishes along a
link $L_v\subset Y$ with multiplicity two. Being a non--zero section of
$TY|_{Y\setminus L_v}$, $v$ determines a trivialization and so a $Spin$
structure on $Y\setminus L_v$. Since $v$ vanishes with multiplicity two along
$L_v$, this $Spin$ structure extends uniquely to a $Spin$ structure on $Y$,
which induces a $Spin^c$ structure $\t_v\in Spin^c(Y)$. The link $L_v$ carries
a natural orientation such that $2\PD([L_v])$ equals the Euler class of
$\xi$. According to \cite{G} the $Spin^c$ structure $\t_\xi$ is given by
\begin{equation}\label{e:spincform}
\t_{\xi}=\t_v+\PD([L_v]).
\end{equation}

\begin{lem}\label{l:chern-section}
Let $n\geq 2g$, and let $\xi_0$ and $\xi_1$ be the contact structures
on $Y_{g,n}$ given by Definition~\ref{d:contact}. Let $F\in
H^2(Y_{g,n};\Z)$ denote the Poincar\'e dual of the homology class of a
fiber of the fibration $\pi\co Y_{g,n}\to\Si_g$. Then, the Euler class
of $\xi_i$, $i\in\{0,1\}$, as an oriented two--plane bundle is equal
to
\[
(-1)^{i+1}2gF.
\]
If $n$ is even, then $\xi_i$ admits a section $v$ vanishing
with multiplicity two along a smooth curve $L_v\subset Y_{g,n}$ with
\[
PD([L_v])= (-1)^{i}({\frac n2} -g)F.
\]
Moreover, the $Spin$ structure $\t_v\in Spin(Y_{g,n})$ is equal to
$\t_{g+\frac{n}{2}-1}$.
\end{lem}

\begin{proof}
The contact structure $\ze$ on $Y_{g,2g-2}$ given by
Lemma~\ref{l:horizontal} is tangent to the fibers of the fibration
$Y_{g,2g-2}\to\Si_g$. Therefore, any nowhere vanishing vector field
tangent to the fibers gives a nowhere zero section $v$ of $\ze$. By
the construction~\eqref{e:contactgluing} defining $\xi_i$, this gives
a nowhere vanishing section $v$ of $\xi_i\vert_{Y_{g,2g-2}\setminus
U}$, which glues up to the section given in (2) of
Lemma~\ref{l:torus-section}. Therefore, by Definition~\ref{d:contact}
the Euler class of $\xi_i$ is 
\[
(-1)^{i}(p-1)F=(-1)^{i}(n-2g)F=(-1)^{i+1}2gF.
\]
If $n$ is even, then $p-1=n-2g$ is even as well. In this case we may
assume that the section given in (2) of Lemma~\ref{l:torus-section}
vanishes with multiplicity two along a smooth curve representing
$(-1)^{i}(\frac {p-1}{2})=(-1)^{i}({\frac n2} -g)$ times the homology
class of the core circle. Splicing such a section to the nowhere zero
section $v$ we obtain the formula for $PD([L_v])$.

We see from Equation~\eqref{e:spincform} that if $v$ is the non--vanishing
section of $\ze$ used above, we have $\t_v=\t_\ze$. On the other hand, by
Remark~\ref{r:transverse} $\t_\ze$ is equal to $\t_0$, which extends as a
$Spin$ structure to $D_{g,n}$ by definition. Thus, the $Spin$ structure
$\t_v\in Spin(Y_{g,n})$ extends to $D_{g,n}$ away from the preimage
$\pi^{-1}(D^2)$ of some two--disk $D^2\subset\Si_g$. But $\pi^{-1}(D^2)$ is
homeomorphic to a ball, therefore the unique $Spin$ structure on $\del
\pi^{-1}(D^2)$ extends for trivial reasons to the unique $Spin$ structure on
$\pi^{-1}(D^2)$.  This proves that $\t_v$ extends to $D_{g,n}$ as a $Spin$
structure. Therefore, by Remark~\ref{r:ext} $\t_v$ must coincide with
$\t_{g+\frac{n}{2}-1}$.
\end{proof} 

\begin{prop}\label{p:spinc}
Let $n\geq 2g$ and let $\xi_0$ and $\xi_1$ be the contact structures
on $Y_{g,n}$ from Definition~\ref{d:contact}. Then,
$\t_{\xi_0}=\t_{n-1}$ and $\t_{\xi_1}=\t_{2g-1}$.
\end{prop}

\proof
By Lemma~\ref{l:chern-section}, the Euler class of $\xi_i$ coincides
with $c_1(\t_{2ig-1})$, $i=0,1$. If $n$ is odd, $H^2(Y_{g,n};\Z)$ has
no $2$--torsion and the result follows. If $n$ is even, by
Lemma~\ref{l:chern-section} and Equation~\eqref{e:spincform} we have
$$
\t_{\xi_i}=\t_{g+\frac{n}{2}-1} + (-1)^{i}({\frac n2} -g)F =
\t_{2ig-1}.
\eqno{\qed}$$

\begin{rem}\label{r:homotopy}
  Since the $Spin^c$ structures $\t_{\xi_i}$ are homotopy invariants,
  Proposition~\ref{p:spinc} implies that $\xi_0$ and $\xi_1$ are not
  homotopic as oriented tangent two--plane fields on $Y_{g,n}$ once
  $n>2g$.  It can be shown that $\xi_0$ and $\xi_1$ are
  contactomorphic, i.e.~there is a self--diffeomorphism of $Y_{g,n}$
  sending one to the other. We will not use that fact in this paper.
\end{rem}

\section{Monopole equations and the proof of Theorem~\ref{t:main}}
\label{s:gauge}

This section is devoted to the proof of the main result of the paper,
Theorem~\ref{t:main}. For definitions and properties of the solutions
to the Seiberg--Witten equations on cylinders $\R\x Y$ and on
symplectic fillings, we refer the reader to~\cite{MOY}
and~\cite{KrMr}.

\begin{lem}\label{l:perturb}
Let $n\geq 2g$, and let $\xi_i$, $i\in\{0,1\}$ be one of the contact
structures on $Y_{g,n}$ from Definition~\ref{d:contact}. Then,
\begin{enumerate}
\item
The moduli space $N(Y_{g,n},\t_{\xi_i})$ of solutions to the
unperturbed Seiberg--Witten equations is smooth and consists of
reducibles.
\item There exists a real number $\ep>0$ such that, if
  $\mu\in\Om^2(Y_{g,n})$ is a closed two--form whose $L^2$--norm is
  smaller than $\ep$, then either 
\[
[\mu]= 2\pi c_1(\t_{\xi_i})\in H^2(Y_{g,n};\R)
\]
or the $\mu$--perturbed Seiberg--Witten moduli space
$N_\mu(Y_{g,n},\t_{\xi_i})$ is empty.
\end{enumerate}
\end{lem}

\begin{proof}
By Theorems~1 and~2 of~\cite{MOY} (see also \cite[Theorem~2.2]{OSz}), if
$|n|>2g-2$ and $e\not\in [0,2g-2]$, then the moduli space $N(Y_{g,n},\t_e)$ is
smooth and contains only reducibles. Therefore, part (1) of the statement
follows from Proposition~\ref{p:spinc}.

To prove (2) we argue by contradiction. Let
$\{\mu_n\} _{n=1}^{\infty}
\subset\Om^2(Y_{g,n})$ be a sequence of closed two--forms
such that $[\mu_n]\neq 2\pi c_1(\t_{\xi_i})$,
$\|\mu_n\|\overset{n\to\infty}{\to} 0$, and $N_{\mu_n}(Y,\t_{\xi_i})$
contains some element $[(A_n,\psi_n)]$. By a standard compactness
argument there is a subsequence converging, modulo gauge, to a
solution $(A_0,\psi_0)$ of the unperturbed Seiberg--Witten equations, 
which is reducible by part (1). We have $\psi_n\neq 0$ because the
assumption $[\mu_n]\neq 2\pi c_1(\t_{\xi_i})$ implies that
$(A_n,\psi_n)$ must be irreducible, so we can set
$\varphi_n=\frac{\psi_n}{\|\psi_n\|}$ for every $n$.

The smoothness of $N(Y_{g,n},\t_{\xi_i})$ implies that the kernel of
the Dirac operator $D_{A_0}$ is trivial, so by standard elliptic
estimates (see e.g.~\cite[page~423]{DK}) we have
\begin{equation}\label{e:ell}
1=\|\varphi_n\|\leq C\| D_{A_0}\varphi_n\|
\end{equation}
for some constant $C$. On the other hand, since 
the $(A_n,\psi_n)$'s are solutions to the Seiberg--Witten equations,
by writing $A_n=A_0+a_n$ we have
\[
0=D_{A_n}\varphi_n = D_{A_0}\varphi_n + a_n\cdot\varphi_n,
\]
where `$\cdot$' denotes Clifford multiplication. Since
$\|\varphi_n\|=1$ while $\|a_n\|\to 0$, this implies
$D_{A_0}\varphi_n\to 0$, contradicting~\eqref{e:ell}.
\end{proof}

\begin{prop}\label{p:filling}
Let $n\geq 2g$, and let $\xi_i$, $i\in\{0,1\}$ be one of the contact
structures on $Y_{g,n}$ from Definition~\ref{d:contact}. Let $(W,\om)$
be a weak symplectic semi--filling of $(Y_{g,n},
\xi_i)$. Then $\del W$ is connected, $b_2 ^+(W)=0$ and the
homomorphism 
\[
H^2 (W;\R )\to H^2(\partial W;\R) 
\]
induced by the inclusion $\del W\subset W$ is the zero map.
\end{prop}

\begin{proof}
In~\cite{PLpos} it is proved that if $(W,\om)$ is a weak semi--filling of a
contact three--manifold $(Y,\xi)$, where $Y$ carries metrics with positive
scalar curvature, then (a) $\del W$ is connected and $b_2 ^+(W)=0$
(\cite[Theorem~1.4]{PLpos}) and (b) the homomorphism 
\[
H^2 (W;\R )\to H^2(\partial W;\R) 
\]
induced by the inclusion $\del W\subset W$ is the zero map
(\cite[Proposition~2.1]{PLpos}).

The positive scalar curvature assumption was used in the proof of
\cite[Theorem~1.4]{PLpos} to guarantee that the moduli space
$N(Y,\t_\xi)$ of solutions to the unperturbed Seiberg--Witten
equations is smooth and consists of reducibles. But this is true for
the contact structures $\xi_0$ and $\xi_1$ on $Y_{g,n}$ by part (1) of
Lemma~\ref{l:perturb}. Therefore, conclusion (a) holds for any
symplectic semi--filling of $(Y_{g,n},\xi_i)$, $i\in\{0,1\}$.

Similarly, the existence of positive scalar curvature metrics was used
in \cite[Proposition~2.1]{PLpos} to prove that if $\mu$ is a closed
two--form whose $L^2$--norm is sufficiently small and $[\mu]\neq 2\pi
c_1(\t_\xi)$, then the $\mu$--perturbed Seiberg--Witten moduli space
$N_\mu(Y,\t_\xi)$ is empty. But for the contact structures $\xi_0$ and
$\xi_1$ this is precisely the content of (2) in Lemma~\ref{l:perturb}.
Hence, conclusion (b) holds for symplectic semi--fillings of
$(Y_{g,n},\xi_i)$, $i\in\{0,1\}$. 
\end{proof}

We shall now give a purely topological argument showing that, under
some restrictions on $g$ and $n$, there is no smooth four--manifold
$W$ with $\del W=Y_{g,n}$ satifying the conclusion of
Proposition~\ref{p:filling}.

Let $\CP^2$ be the complex projective plane. Denote by $\widehat{\CP^2}$
the blow-up of $\CP^2$ at $k$ distinct points. The second homology
group $H_2(\widehat{\CP^2};\Z)$ is generated by classes
$h,e_1,e_2,\ldots,e_k$, where $h$ corresponds to the standard
generator of $H_2(\CP^2;\Z)$ and the $e_i$'s are the classes of the
exceptional curves.

Let $d$ be a positive integer, and suppose $k\geq 2d$. Define
$\L_d=(H_d,Q_d)$ as the intersection lattice given by the subgroup
\[
H_d=\langle e_1-e_2, e_2-e_3, \ldots , e_{2d-1}-e_{2d},
h-e_1-e_2-\ldots - e_d\rangle\subset H_2(\widehat{\CP^2};\Z)
\] 
together with the restriction $Q_d$ of the intersection form
$Q_{\widehat{\CP^2}}$. For $m\geq 1$ let $\D_m=(\Z^m, m(-1))$ be
that standard negative definite diagonal lattice. The following lemma
generalizes a result from~\cite{Paolo}.

\begin{lem}\label{l:nonemb}
$\L_d$ does not embed into $\D_m$.
\end{lem}
\begin{proof} 
Arguing by contradiction, suppose that $j\co\L_d\hookrightarrow\D_m$
is an embedding. Then, 
\[
\rk \D_m=m\geq \rk\L_d=2d. 
\]
Let $w_1,\ldots ,w_{2d}$ be the obvious generators of $\L_d$ with
$w_i\cdot w_i=-2$ for $i\leq 2d-1$ and $w_{2d}\cdot
w_{2d}=-(d-1)$. Then, since the elements $w_1,\ldots, w_{2d-1}$ have
square $-2$, there is a set $e_1,\ldots,e_m$ of standard generators of
$\D_m$ such that 
\[
j(w_i)=e_i-e_{i+1},\quad i=1,\ldots,2d-1. 
\]
Suppose that 
\[
j(w_{2d})=a_1e_1+\cdots +a_m e_m,\quad a_1,\ldots, a_m\in\Z.
\]
Then, the relations $w_{2d}\cdot w_i =0$ for $i\neq d$ and
$w_{2d}\cdot w_d=1$ imply
\[
j(w_{2d})=a\sum_{i=1}^{2d} e_i + \sum_{i=d+1}^{2d} e_i + \sum_{i\geq
2d+1} a_i e_i,
\]
for some $a\in\Z$, therefore $j(w_{2d})\cdot j(w_{2d})\leq -d$, which
is impossible because $w_{2d}\cdot w_{2d}=-(d-1)$.
\end{proof}

\begin{prop}\label{p:nonexist}
Suppose that $2g\geq d(d+1)$ and $n\leq (d+1)^2+3$ for some positive
integer $d$. Then, there is no smooth, compact, oriented
four--manifold $W$ such that $\del W=Y_{g,n}$, $b_2^+(W)=0$ and the
map 
\[
H^2(W;\R)\to H^2(\del W;\R) 
\]
is the zero map.
\end{prop}

\begin{proof}
Consider a smooth curve $C\subset\CP^2$ of degree $d+2$, and let
$\widehat{\CP^2}$ be the blow-up of $\CP^2$ at $k$ distinct points of
$C$. Let $\widehat C$ be the proper trasform of $C$ inside
$\widehat{\CP^2}$.

Denote by ${\widetilde C}\subset\widehat{\CP^2}$ a smooth, oriented
surface obtained by adding $g - \frac 12 d(d+1)$ fake handles to
$\widehat C$. Let $Z\subset\widehat{\CP^2}$ be the complement of an
open tubular neighborhood of $\widetilde C$. The boundary of $Z$
is orientation--reversing diffeomorphic to $Y_{g,n}$, where
$n=(d+2)^2-k$. Clearly, any $n\leq (d+1)^2+3$ can be realized by some
$k\geq 2d$, in which case the lattice $\L_d$ can be realized as a
sublattice of the intersection lattice of $\widehat{\CP^2}$. The
generators $e_i-e_{i+1}$ and $h-e_1-\ldots - e_d$ of $\L_d$ can all be
represented by smooth surfaces inside $Z$, so we have an embedding
\[
\L_d\subset(H_2(Z;\Z),Q_Z).
\]
Now we argue by contradiction. Let $W$ be a smooth four--manifold as
in the statement. Consider a smooth, closed four--manifold $V$ of
the form 
\[
W\cup _{Y_{g,n}} Z. 
\]
Since 
\[
H^2(W;\R)\to H^2(\del W;\R)
\]
is zero, by Poincar\'e duality the map 
\[
H_2(W,\del W;\R)\to H_1(\del W;\R)
\]
vanishes. This implies that 
\[
H_1(\del W;\R)\to H_1(W;\R) 
\]
is injective and by Mayer--Vietoris 
\[
H_2(V;\R)=H_2(W;\R) + H_2(Z;\R).
\]
Therefore, since $b_2^+(W)=b_2^+(Z)=0$, we have $b_2^+(V)=0$.
Donaldson's theorem on the intersection form of closed, definite
four--manifolds (see~\cite{Do1, Do2}) implies that the intersection
lattice $(H_2(V;\Z),Q_V)$ must be isomorphic to $\D_{b_2(V)}$. Thus,
the resulting existence of an embedding 
\[
\L_d\subset(H_2(Z;\Z),Q_Z)
\]
contradicts Lemma~\ref{l:nonemb}.
\end{proof}

\begin{proof}[Proof of Theorem~\ref{t:main}]
The statement follows from Propositions~\ref{p:filling}
and~\ref{p:nonexist}.
\end{proof}

\begin{rem}
The restrictions on $n$ and $g$ appearing in the statement of
Proposition~\ref{p:nonexist} are only due to the inability of the
authors to construct more smooth four--manifolds with the necessary
properties.  In fact, using slightly different methods, in \cite{LS}
we proved that the conclusion of Theorem~\ref{t:main} holds for every
$n \geq 2g>0$.
\end{rem}

\end{document}